\date{November 4, 2006}

\documentclass[11pt,leqno]{article}
\usepackage{graphicx}
\usepackage{amssymb}
\usepackage{amsmath}
\newtheorem{definition}{Definition}[section]
\newtheorem{lemma}[definition]{Lemma}
\newtheorem{theorem}[definition]{Theorem}
\newtheorem{proposition}[definition]{Proposition}

\newtheorem{remark}[definition]{Remark}

 1

\newenvironment{proof}{\noindent{\bf Proof~:}}{\QED\medskip}
\def\QED{\hskip0.1em\hfill\null\ \null\nobreak\hfill
\kern3pt\lower1.8pt\vbox{\hrule\hbox
{\vrule\kern1pt\vbox{\kern1.7pt \hbox{$\scriptstyle
QED$}\kern0.2pt}\kern1pt\vrule}\hrule}}

\renewcommand{\a}{\alpha}
\renewcommand{\b}{\beta}

\newcommand{\bmu}{\bar{\mu}}
\newcommand{\bnu}{\bar{\nu}}
\newcommand{\baeta}{\bar{\eta}}
\newcommand{\btheta}{\bar{\theta}}

\newcommand{\PP}{{\mathbb{P}}}
\newcommand{\ZZ}{{\mathbb{Z}}}
\newcommand{\CC}{{\mathbb{C}}}

\newcommand{\cO}{{\mathcal{O}}}

\newcommand{\iso}{\cong}

\newcommand{\surj}{\twoheadrightarrow}
\newcommand{\inc}{\hookrightarrow}

\newcommand{\la}{\langle}
\newcommand{\ra}{\rangle}

\begin{document}

\title{An $8$-dimensional non-formal simply connected symplectic manifold}

\author{Marisa Fern\'andez and Vicente Mu\~noz}

\maketitle

\begin{abstract}
 We answer in the affirmative the question posed by Babenko
 and Taimanov \cite{BT2} on the existence of
 non-formal simply connected compact symplectic manifolds
 of dimension $8$.
\end{abstract}

\section{Introduction} \label{sec:introduction}

Simply connected compact manifolds of dimension less than or equal
to $6$ are formal \cite{NM,FM1}, and there are simply connected
compact manifolds of dimension greater than or equal to $7$ which
are non-formal \cite{Oprea,FM2,Rudyak,Cav1,FM3}. If we are
treating the symplectic case, the story is not so straightforward.
Lupton and Oprea \cite{LO} conjectured that any simply connected
compact symplectic manifold is formal. Babenko and Taimanov
\cite{BT1,BT2} disproved this conjecture giving examples  of
non-formal simply connected compact {\em symplectic} manifolds of
any dimension bigger than or equal to $10$, by using the
symplectic blow-up \cite{McDuff}. They raise the question of the
existence of non-formal simply connected compact symplectic
manifolds of dimension $8$. The techniques of construction of
symplectic manifolds used so far
\cite{BT,BT2,Cav1,Cav2,FM1,IRTU,RT,TO} have not proved fruitful when
addressing this problem. In this note, we answer the question in
the affirmative by proving the following.

\begin{theorem} \label{thm:main}
 There is a simply connected compact symplectic manifold
 of dimension $8$ which is non-formal.
\end{theorem}

To construct such a manifold, we introduce a new technique to
produce symplectic manifolds, which we hope can be useful for
obtaining examples with interesting properties. We consider a
non-formal compact symplectic $8$-dimensional manifold with a
symplectic non-free action of a finite group such that the
quotient space is a non-formal orbifold which is simply connected.
Then we resolve symplectically the singularities to produce a {\em
smooth\/} symplectic $8$-manifold satisfying the required
properties. The origin of the idea stems from our study of Guan's
examples \cite{Gu} of compact holomorphic symplectic manifolds
which are not K\"ahler.

\section{A simply-connected symplectic $8$-manifold} \label{sec:manifold}

Consider the complex Heisenberg group $H_{\CC}$, that is, the
complex nilpotent Lie group of complex  matrices of the form
 $$
 \begin{pmatrix} 1&u_2&u_3\\ 0&1&u_1\\ 0&0&1\end{pmatrix},
 $$
and let $G=H_\CC \times \CC$, where $\CC$ is the additive group of
complex numbers. We denote by $u_4$ the coordinate function
corresponding to this extra factor. In terms of the natural
(complex) coordinate functions $(u_1,u_2,u_3,u_4)$ on $G$, we have
that the complex $1$-forms $\mu=du_1$, $\nu=du_2$,
$\theta=du_3-u_2 \, du_1$ and $\eta=du_4$ are left invariant, and
 $$
 d\mu=d\nu=d\eta=0, \quad d\theta=\mu\wedge\nu.
 $$

Let $\Lambda \subset \CC$ be the lattice generated by $1$
and $\zeta=e^{2\pi i/3}$,
and consider the discrete subgroup
$\Gamma\subset G$ formed by the matrices in which $u_1,u_2,u_3,
u_4 \in \Lambda$. We define the compact (parallelizable)
nilmanifold
 $$
 M=\Gamma \backslash G.
 $$
We can describe $M$ as a principal torus bundle
 $$
 T^2=\CC/\Lambda \inc M \to T^6=(\CC/\Lambda)^3,
 $$
by the projection $(u_1,u_2,u_3,u_4) \mapsto (u_1,u_2,u_4)$.

Now introduce the following action of the finite group $\ZZ_{3}$
 \begin{eqnarray*}
 \rho: G & \to & G\\
 (u_1, u_2,u_3, u_4) &\mapsto & (\zeta\, u_1, \zeta\, u_2,\zeta^2\, u_3,\zeta\, u_4) .
 \end{eqnarray*}
This action satisfies that $\rho(p\cdot q)=\rho(p)\cdot \rho(q)$,
for $p, q \in G$, where the dot denotes the natural group
structure of $G$. The map $\rho$ is a particular case of a
homothetic transformation (by $\zeta$ in this case) which is well
defined for all nilpotent simply connected Lie groups with graded
Lie algebra. Moreover $\rho(\Gamma)=\Gamma$, therefore $\rho$
induces an action on the quotient $M=\Gamma\backslash G$.  The
action on the forms is given by
 $$
 \rho^*\mu= \zeta\, \mu, \quad
 \rho^*\nu= \zeta\, \nu,\quad
 \rho^*\theta= \zeta^2\,\theta, \quad
 \rho^*\eta= \zeta\,\eta.
 $$

The complex $2$-form
 $$
 \omega= i \,\mu \wedge \bmu + \nu\wedge \theta +
 \bnu \wedge \btheta+ i \,\eta \wedge \baeta
 $$
is actually a real form which is clearly closed and
which satisfies $\omega^4\not=0$. Thus $\omega$ is a symplectic
form on $M$. Moreover, $\omega$ is $\ZZ_3$-invariant. Hence the
space
 $$
 \widehat M=M/\ZZ_{3}
 $$
is a symplectic orbifold, with the symplectic form
$\widehat\omega$ induced by $\omega$. Our next step is to find a
smooth symplectic manifold $\widetilde M$ that desingularises
$\widehat M$.

\begin{proposition}\label{prop:desingularization}
 There exists a smooth compact
 symplectic manifold $(\widetilde{M},\widetilde{\omega})$
 which is isomorphic to
 $(\widehat{M},\widehat{\omega})$ outside the singular points.
\end{proposition}

\begin{proof}
Let $p \in M$ be a fixed point of the $\ZZ_3$-action. Translating
by a group element $g\in G$ taking $p$ to the origin, we may
suppose that $p=(0,0,0,0)$ in our coordinates. At $p$, the
symplectic form is
 $$
  \omega_0=i \,du_1\wedge d\bar{u}_1+du_2 \wedge du_3+
  d\bar{u}_2\wedge d\bar{u}_3+i\, du_4\wedge d\bar{u}_4.
 $$
Take now $\ZZ_3$-equivariant Darboux coordinates around $p$, $\Phi
\colon (B,\omega) \longrightarrow
(B_{\CC^4}(0,\epsilon),\omega_0)$, for some $\epsilon>0$. This
means that $\Phi^*\omega_0=\omega$ and $\Phi\circ
\rho=d\rho_{p}\circ \Phi$, where we interpret
$(B_{\CC^4}(0,\epsilon),\omega_0)\subset (T_pM,\omega_0)\cong
(\CC^4,\omega_0)$ in the natural way. (The proof of the existence
of usual Darboux coordinates in \cite[pp.\ 91--93]{McDuff-Salamon}
carry over to this case, only being careful that all the objects
constructed should be $\ZZ_3$-equivariant.) We denote the new
coordinates given by $\Phi$ as $(u_1,u_2,u_3,u_4)$ again (although
they are not the same coordinates as before).

Now introduce the new set of coordinates:
 $$
 (w_1,w_2,w_3,w_4)=
 (u_1,\frac{1}{\sqrt2}(u_2+\bar{u}_3),
 \frac{i}{\sqrt2}(u_3-\bar{u}_2),u_4).
 $$
Then the symplectic form $\omega$ can be expressed as
  $$
   \omega= i\, (dw_1\wedge d\bar{w}_1+dw_2\wedge d\bar{w}_2
   +dw_3\wedge d\bar{w}_3+  dw_4\wedge d\bar{w}_4 ).
  $$
Moreover, with respect to these coordinates, the $\ZZ_3$-action
$\rho$ is given as
  $$
   \rho(w_1,w_2,w_3,w_4)=(\zeta\, w_1,\zeta\, w_2,\zeta^2\, w_3,\zeta\,
   w_4).
  $$

With this K\"ahler model for a neighbourhood $B$ of $p$, we may
resolve the singularity of $B/\ZZ_3$ with a non-singular K\"ahler
model. Basically, blow up $B$ at $p$ to get $\widetilde{B}$. This
replaces the point with a complex projective space $\PP^3$ in
which $\ZZ_3$ acts as
  $$
   [w_1,w_2,w_3,w_4] \mapsto
   [\zeta\, w_1,\zeta\, w_2,\zeta^2\, w_3,\zeta\, w_4]=
   [w_1,w_2,\zeta\, w_3,w_4].
  $$
Therefore there are two components of the fix-point locus of the
$\ZZ_3$-action on $\widetilde{B}$, namely the point $q=[0,0,1,0]$
and the complex projective plane $H=\{[w_1,w_2,0,w_4]\} \subset
F=\PP^3$. Next blow up $\widetilde{B}$ at $q$ and at $H$ to get
$\widetilde{\widetilde{B}}$. The point $q$ is substituted by a
projective space $H_1=\PP^3$. The normal bundle of $H\subset
\widetilde{B}$ is the sum of the normal bundle of $H\subset F$,
which is $\cO_{\PP^2}(1)$, and the restriction of the normal
bundle of $F\subset\widetilde{B}$ to $H$, which is
$\cO_{\PP^3}(-1)|_{\PP^2}= \cO_{\PP^2}(-1)$. Therefore the second
blow-up replaces the plane $H$ by the $\PP^1$-bundle over $\PP^2$
defined as $H_2 = \PP( \cO_{\PP^2}(-1)\oplus\cO_{\PP^2}(1))$. The
strict transform of $F\subset\widetilde{B}$ under the second
blow-up is the blow up $\widetilde{F}$ of $F=\PP^3$ at $q$, which
is a $\PP^1$-bundle over $\PP^2$, actually $\widetilde{F} = \PP(
\cO_{\PP^2}\oplus\cO_{\PP^2}(1))$. See Figure 1 below.

\begin{figure}[ht]
\begin{center}
\includegraphics{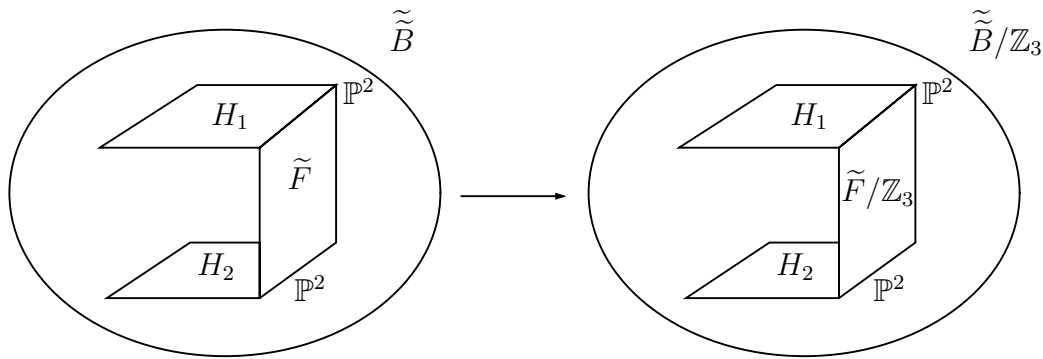}
\end{center}
\caption{Desingularisation process}
\end{figure}

The fix-point locus of the $\ZZ_3$-action on
$\widetilde{\widetilde{B}}$ are exactly the two disjoint divisors
$H_1$ and $H_2$. Therefore the quotient
$\widetilde{\widetilde{B}}/\ZZ^3$ is a smooth K\"ahler manifold
\cite[page 82]{BPV}. This provides a symplectic resolution of the
singularity $B/\ZZ_3$. To glue this K\"ahler model to the
symplectic form in the complement of the singular point we use
Lemma \ref{lem:glue} below. We do this at every fixed point to get
a smooth symplectic resolution of $\widehat{M}$.
\end{proof}

\begin{lemma}\label{lem:glue}
 Let $(B,\omega_0)$ be the standard K{\"a}hler ball in $\CC^n$, $n>1$, and let $\Pi$
 be a finite group acting linearly (by complex isometries)
 on $B$ whose only fixed point is the origin.
 Let $\phi:(\widetilde{B},\omega_1) \to (B/\Pi,\omega_0)$ be a K{\"a}hler
 resolution of the singularity of the quotient. Then there is a
 symplectic form $\Omega$ on $\widetilde{B}$ which coincides with
 $\omega_0$ near the boundary, and with a positive multiple of
 $\omega_1$ near the exceptional divisor $E=\phi^{-1}(0)$.
 Moreover $\Omega$ is tamed by the complex structure.
\end{lemma}

\begin{proof}
Since $\phi:(\widetilde{B},\omega_1) \to (B/\Pi,\omega_0)$ is
holomorphic, $\omega_0$ and $\omega_1$ are K{\"a}hler forms in
$\widetilde{B} - E = B-\{0\}$ with respect to the same complex
structure $J$. Therefore $(1-t)\omega_0+ t\omega_1$ is a K{\"a}hler
form on $\widetilde{B}$, for any number $0<t<1$. (Note that
$\omega_0|_E=0$, where we denote again by $\omega_0$ the pull-back
to $\widetilde{B}$.)

Fix $\delta>0$ small and let $A=\{z\in B\, |\, \delta
<|z|<2\delta\} \subset B$. Since $A$ is simply connected, we may
write $\omega_1-\omega_0=d\alpha$, with $\alpha\in \Omega^1(A)$,
which we can furthermore suppose $\Pi$-invariant.

Let $\rho: [0,\infty)\to [0,1]$ be a smooth function whose value
is $1$ for $r\leq 1.1\delta$ and $0$ for $r\geq 1.9\delta$. Define
  $$
  \Omega=\omega_0 + \epsilon \ d(\rho (|z|) \alpha).
  $$
This equals $\omega_0$ for $|z|\geq 1.9\delta$, and
$\omega_0+\epsilon (\omega_1-\omega_0)=(1-\epsilon) \omega_0
+\epsilon \ \omega_1$ for $|z|\leq 1.1\delta$. For $1.1\delta \leq
|z| \leq 1.9\delta$, let $C>0$ be a bound of $d(\rho \alpha)
(u,Ju)$, for any $u$ unitary tangent vector (with respect to the
K{\"a}hler form $\omega_0$). Choose $0<\epsilon<\min\{1,C^{-1}\}$.
Then $\Omega(u,Ju)>0$ for any non-zero $u$.
\end{proof}

\begin{proposition}
 The manifold $\widetilde{M}$ is simply connected.
\end{proposition}

\begin{proof}
Fix the base points: let $p_0 \in M=\Gamma\backslash G$ be the
image of $(0,0,0,0)\in G$ and let $\hat{p}_0\in \widehat M$ be the
image of $p_0$ under the projection $M\to \widehat M$. There is an
epimorphism of fundamental groups
 $$
 \Gamma=\pi_1(M) \surj \pi_1(\widehat M),
 $$
since the $\ZZ_3$-action has a fixed point \cite[Corollary
6.3]{Bredon}. Now the nilmanifold $M$ is a principal $2$-torus
bundle over the $6$-torus $T^6$, so we have an exact sequence
 $$
 \ZZ^2 \inc \Gamma \to \ZZ^6.
 $$
Let $\bar{p}_0=\pi(p_0)$, where $\pi \colon M \to T^6$ denotes the
projection of the torus bundle. Clearly, $\ZZ_3$ acts on
$\pi^{-1}(\bar{p}_0)\iso T^2=\CC/\Lambda$ with  $3$ fixed points,
and the quotient space $T^2/\ZZ_{3}$ is a $2$-sphere $S^2$. So the
restriction to $\ZZ^2=\pi_1(T^2)\subset \pi_1(M)=\Gamma$ of the
map $\Gamma \surj \pi_1(\widehat M)$ factors through
$\pi_1(T^2/{\ZZ_3})=\{1\}$, hence it is trivial. Thus the map
$\Gamma \surj \pi_1(\widehat M)$ factors through the quotient
$\ZZ^6 \surj \pi_1(\widehat M)$. But $M$ contains three
$\ZZ_3$-invariant $2$-tori, $T_1$, $T_2$ and $T_3$ (which are the
images of $\{ (u_1,0,0,0)\}$, $\{ (0,u_2,0,0)\}$ and $\{
(0,0,0,u_4)\}$, respectively) such that $\pi_1(\widehat M)$ is
generated by the images of $\pi_1(T_1)$, $\pi_1(T_2)$ and
$\pi_1(T_3)$. Again, each quotient $T_i/\ZZ_3$ is a $2$-sphere,
hence $\pi_1(\widehat M)$ is generated by
$\pi_1(T_i/\ZZ_3)=\{1\}$, which proves that $\pi_1(\widehat
M)=\{1\}$.

Finally, the resolution $\widetilde M\to \widehat M$ consists of
substituting, for each singular point $p$, a neighbourhood
$B/\ZZ_3$ of it by a non-singular model $\widetilde{\widetilde
B}/\ZZ_3$. The fiber over the origin of $\widetilde{\widetilde
B}/\ZZ_3\to B/\ZZ_3$ is simply connected: it consists of the union
of the three divisors $H_1=\PP^3$, $H_2=\PP(
\cO_{\PP^2}(-1)\oplus\cO_{\PP^2}(1))$ and $\widetilde F/\ZZ_3=\PP(
\cO_{\PP^2}\oplus\cO_{\PP^2}(3))$, all of them are simply
connected spaces, and their intersection pattern forms no cycles
(see Figure 1). Therefore, a simple Seifert-Van Kampen argument
proves that $\widetilde{M}$ is simply connected
\end{proof}

\begin{lemma}\label{lem:cohomology}
The odd degree Betti numbers of $\widetilde M$ are $b_1(\widetilde
M)=b_3(\widetilde M)=b_5(\widetilde M)=b_7(\widetilde M)=0$.
\end{lemma}

\begin{proof}
As $\widetilde M$ is simply connected, then $b_1(\widetilde M)=0$.
Next, using Nomizu's theorem \cite{No} to compute the cohomology
of the nilmanifold $M$, we easily find that $H^3(M)=W\oplus
\overline{W}$, where
 \begin{align*}
 W= \la &[\mu\wedge\bmu\wedge \eta],  [\nu\wedge\bnu\wedge
 \eta], [\mu\wedge\bnu\wedge\eta], [\bmu \wedge \nu\wedge\eta],
 [\mu\wedge\eta\wedge\baeta],[\nu\wedge\eta\wedge\baeta],
  [\mu\wedge\nu\wedge\theta],\\ &[\mu\wedge\bnu\wedge\btheta],[\bmu\wedge\nu\wedge\btheta],
 [\mu\wedge\bmu\wedge\btheta],[\nu\wedge\bnu\wedge\btheta],
 [\mu\wedge\eta\wedge\theta],[\nu\wedge\eta\wedge\theta],[\bmu\wedge\eta\wedge\btheta],
 [\bnu\wedge\eta\wedge\btheta] \ra
 \end{align*}
and $\overline W$ is its complex conjugate. (Here $H^*(X)$ denotes
cohomology with complex coefficients.) Clearly $\rho$ acts as
multiplication by $\zeta$ on $W$ and as multiplication by
$\zeta^2=\bar\zeta$ on $\overline W$. Therefore
$H^3(\widehat{M})=H^3(M)^{\ZZ_3}=0$.

The desingularisation process of Proposition
\ref{prop:desingularization} consists on removing contractible
neighborhoods of the form $B_i/\ZZ_3$, $B_i\cong
B_{\CC^4}(0,\epsilon)$, around each fixed point $p_i$, and
inserting a non-singular K\"ahler model
$\widetilde{\widetilde{B}}_i/\ZZ_3$ which retracts to the
``exceptional divisor'' $E_i=\phi^{-1}(0)$,
$\phi:\widetilde{\widetilde{B}}_i/\ZZ_3\to B_i/\ZZ_3$. We glue
along the region $A/\ZZ_3$ which retracts into $S^7/\ZZ_3$, a
rational homology $7$-sphere. An easy Mayer-Vietoris argument then
shows that $H^j(\widetilde{M})=H^j(\widehat{M})\oplus (\bigoplus_i
H^j(E_i))$ for $0<j<7$. All the $E_i$ are diffeomorphic to the
$6$-dimensional complex manifold depicted in Figure 1, which
consists of the union of $H_1=\PP^3$, $H_2=\PP(
\cO_{\PP^2}(-1)\oplus\cO_{\PP^2}(1))$ (a $\PP^1$-bundle over
$\PP^2$) and $\widetilde F/\ZZ_3=\PP(
\cO_{\PP^2}\oplus\cO_{\PP^2}(3))$ (another $\PP^1$-bundle over
$\PP^2$), intersecting in copies of the complex projective plane.
So $H^3(E_i)=0$ and hence $H^3(\widetilde{M})=0$.

The statement $b_5(\widetilde M)=b_7(\widetilde M)=0$ follows from
Poincar{\'e} duality.
\end{proof}

\section{Non-formality of the constructed manifold} \label{sec:formality}

Formality for a simply connected manifold $M$ means that its
rational homotopy type is determined by its cohomology algebra.
Let us recall its definition (see \cite{DGMS,TO} for more
details). Let $X$ be a simply connected smooth manifold and
consider its algebra of differential forms $(\Omega^*(X),d)$. Let
$\psi:(\bigwedge V,d)\to (\Omega^*(X),d)$ be a minimal model for
this algebra \cite{DGMS}. 
Then $X$ is formal if there is a quasi-isomorphism $\psi':(\bigwedge
V,d) \to (H^*(X),d=0)$, i.e. a morphism of differential algebras,
inducing the identity on cohomology.

\begin{lemma}\label{lem:nonformal}
Let $X$ be a simply connected smooth manifold with $H^3(X)=0$, and
let $a, x_1,x_2,x_3 \in H^2(X)$ be cohomology classes satisfying
that $a\cup x_i=0$, $i=1,2,3$. Choose forms $\alpha,\beta_i\in
\Omega^2(X)$ and $\xi_i\in \Omega^3(X)$, with $a=[\alpha]$,
$x_i=[\beta_i]$ and $\alpha\wedge \beta_i= d\xi_i$, $i=1,2,3$. If
the cohomology class
 \begin{equation}\label{eqn:GM2}
 [ \xi_1\wedge\xi_2\wedge\beta_3
  +\xi_2\wedge\xi_3\wedge\beta_1
  +\xi_3\wedge\xi_1\wedge\beta_2]\in H^8(X)
 \end{equation}
is non-zero, then $X$ is non-formal.
\end{lemma}

\begin{proof}
First, notice that
 $$
  d(\xi_1\wedge\xi_2\wedge\beta_3
  +\xi_2\wedge\xi_3\wedge\beta_1
  +\xi_3\wedge\xi_1\wedge\beta_2)=\alpha\wedge
 \beta_1\wedge \xi_2\wedge \beta_3 - \xi_1 \wedge \alpha\wedge
 \beta_2 \wedge \beta_3 +
 $$
 $$
 +\alpha\wedge \beta_2\wedge
 \xi_3\wedge \beta_1 - \xi_2 \wedge \alpha\wedge \beta_3 \wedge
 \beta_1 +\alpha\wedge \beta_3\wedge \xi_1\wedge \beta_2 - \xi_3
 \wedge \alpha\wedge \beta_1 \wedge \beta_2 =0 \, ,
 $$
so \eqref{eqn:GM2} is a well-defined cohomology class.

Second, let us see that the cohomology class  \eqref{eqn:GM2} does
not depend on the particular forms $\alpha,\beta_i\in \Omega^2(X)$
and $\xi_i\in \Omega^3(X)$ chosen. If we write $a=[\alpha+d f]$,
with $f\in \Omega^1(X)$, then $(\alpha+df)\wedge \beta_i= d(\xi_i+
f\wedge \beta_i)$ and
  $$
  (\xi_1+f\wedge \beta_1)\wedge(\xi_2+f\wedge \beta_2)\wedge\beta_3
  +(\xi_2+f\wedge \beta_2)\wedge(\xi_3+f\wedge \beta_3)\wedge\beta_1
  +(\xi_3+f\wedge \beta_3)\wedge(\xi_1+f\wedge \beta_1)\wedge\beta_2=
  $$
  $$
  =\xi_1\wedge\xi_2\wedge\beta_3
  +\xi_2\wedge\xi_3\wedge\beta_1
  +\xi_3\wedge\xi_1\wedge\beta_2,
  $$
so the cohomology class (\ref{eqn:GM2}) does not change by
changing the representative of $a$. If we change the
representatives of $x_i$, say for instance $x_1=[\beta_1+ df]$,
$f\in \Omega^1(X)$, then $\alpha\wedge (\beta_1+df)= d(\xi_1+
\alpha\wedge f)$ and
 $$
  (\xi_1+ \alpha \wedge f)\wedge\xi_2\wedge\beta_3
  +\xi_2\wedge\xi_3\wedge(\beta_1+df)
  +\xi_3\wedge(\xi_1 +\alpha\wedge f)\wedge\beta_2=
  $$
  $$
  =\xi_1\wedge\xi_2\wedge\beta_3
  +\xi_2\wedge\xi_3\wedge\beta_1
  +\xi_3\wedge\xi_1\wedge\beta_2 + d(f\wedge\xi_2\wedge \xi_3),
  $$
so the cohomology class (\ref{eqn:GM2}) does not change again.
Finally, if we change the form $\xi_1$ to $\xi_1+ g$, $g\in
\Omega^3(X)$ closed, then
  $$
  (\xi_1+g)\wedge\xi_2\wedge\beta_3
  +\xi_2\wedge\xi_3\wedge\beta_1
  +\xi_3\wedge(\xi_1+g)\wedge\beta_2
  $$
  $$
  =\xi_1\wedge\xi_2\wedge\beta_3
  +\xi_2\wedge\xi_3\wedge\beta_1
  +\xi_3\wedge\xi_1\wedge\beta_2 + g \wedge
  ( \xi_2\wedge\beta_3-\xi_3 \wedge\beta_2),
  $$
and $\xi_2\wedge\beta_3-\xi_3\wedge\beta_2\in \Omega^3(X)$ is
closed, hence exact since $H^3(X)=0$. Also in this case the
cohomology class (\ref{eqn:GM2}) does not change.

To see that $X$ is non-formal, consider the minimal model
$\psi:(\bigwedge V,d) \to (\Omega^*(X),d)$ for $X$. Then there are
closed elements $\hat{a}, \hat{x}_i \in (\bigwedge V)^2$ whose
images are $2$-forms $\alpha, \beta_i$ representing $a,x_i$. Since
$[\hat{a}\cdot \hat{x}_i]=0$, there are elements $\hat{\xi}_i\in
(\bigwedge V)^3$ such that $d\hat{\xi}_i= \hat{a}\cdot \hat{x}_i$.
Let $\xi_i=\psi(\hat{\xi}_i)\in \Omega^3(X)$. So
$d\xi_i=\alpha\wedge\beta_i$, $i=1,2,3$.

If $X$ is formal, then there exists a quasi-isomorphism
$\psi':(\bigwedge V,d) \to (H^*(X),0)$. Note that 
$\psi'(\hat{\xi}_i)=0$ since $H^3(X)=0$. Then
  $$
 [\xi_1\wedge\xi_2\wedge\beta_3
  +\xi_2\wedge\xi_3\wedge\beta_1
  +\xi_3\wedge\xi_1\wedge\beta_2]=\psi'(\hat\xi_1\wedge\hat\xi_2\wedge\hat{x}_3
  +\hat\xi_2\wedge\hat\xi_3\wedge\hat{x}_1
  +\hat\xi_3\wedge\hat\xi_1\wedge\hat{x}_2)=0,
  $$
contradicting our assumption. This proves that $X$ is
non-formal.
\end{proof}

\begin{theorem}\label{thm:nonformal}
 The manifold  $\widetilde M$  is non-formal.
\end{theorem}

\begin{proof}
We start by considering the nilmanifold $M$. Consider the closed
forms:
 $$
 \a = \mu \wedge\bmu, \quad
 \b_1 = \nu\wedge\bnu,\quad
 \b_2 = \nu\wedge\baeta,\quad
 \b_3 = \bnu\wedge\eta .
 $$
Then
 $$
 \a \wedge \b_1 = d(-\theta\wedge\bmu\wedge\bnu),
\quad
 \a \wedge \b_2 = d(-\theta\wedge\bmu\wedge\baeta),
\quad
 \a \wedge \b_3 = d(\btheta\wedge \mu\wedge\eta).
 $$
All the forms $\a$, $\b_1$, $\b_2$, $\b_3$,
$\xi_1=-\theta\wedge\bmu\wedge\bnu$,
$\xi_2=-\theta\wedge\bmu\wedge\baeta$ and $\xi_3=\btheta\wedge
\mu\wedge\eta$ are $\ZZ_3$-invariant. Hence they descend to the
quotient $\widehat M=M/\ZZ_3$. Let $q:M\to \widehat M$ denote the
projection, and define $\hat\a=q_*\a$, $\hat\b_i=q_*\b_i$,
$\hat\xi_i=q_*\xi_i$, $i=1,2,3$. Now take a $\ZZ_3$-equivariant
map $\varphi:M\to M$ which is the identity outside some small
balls around the fixed points, and contracts some smaller balls
into the fixed points. It induces a map
$\hat\varphi:\widehat{M}\to\widehat{M}$ such that
$\hat\varphi\circ q=q\circ \varphi$. The forms
$\tilde\a=\hat\varphi^*\hat\a$,
$\tilde\b_i=\hat\varphi^*\hat\b_i$,
$\tilde\xi_i=\hat\varphi^*\hat\xi_i$, $i=1,2,3$, are zero in a
neighbourhood of the fixed points, therefore they define forms on
$\widetilde M$, by extending them by zero along the exceptional
divisors. Note that $\tilde\a, \tilde\b_i\in
\Omega^2(\widetilde{M})$ are closed forms and $\tilde\xi_i\in
\Omega^3(\widetilde{M})$ satisfies $d\tilde\xi_i=
\tilde\a\wedge\tilde\b_i$, $i=1,2,3$.

By Lemma \ref{lem:cohomology}, $H^3(\widetilde{M})=0$, so we may
apply Lemma \ref{lem:nonformal} to the cohomology classes
$a=[\alpha], b_i=[\beta_i]\in H^2(\widetilde{M})$, $i=1,2,3$. The
cohomology class
 \begin{align*}
 [ \tilde\xi_1\wedge\tilde\xi_2\wedge\tilde\beta_3
  +\tilde\xi_2\wedge\tilde\xi_3\wedge\tilde\beta_1
  +\tilde\xi_3\wedge\tilde\xi_1\wedge\tilde\beta_2] =& \,
  [\hat\varphi^*q_*(\xi_1\wedge\xi_2\wedge\beta_3
  +\xi_2\wedge\xi_3\wedge\beta_1
  +\xi_3\wedge\xi_1\wedge\beta_2)] \\ =& \,
  \hat\varphi^*q_*(2[\theta\wedge
 \mu\wedge\nu\wedge\eta\wedge \btheta\wedge
 \bmu\wedge\bnu\wedge\baeta]) \neq 0,
 \end{align*}
since its integral is
 \begin{align*}
  \int_{\widetilde{M}} \hat\varphi^*q_*(2[\theta\wedge
  \mu\wedge\nu\wedge\eta\wedge \btheta\wedge
  \bmu\wedge\bnu\wedge\baeta])
  &= \int_{\widehat{M}}
  \hat\varphi^*q_*(2[\theta\wedge \mu\wedge\nu\wedge\eta\wedge
  \btheta\wedge \bmu\wedge\bnu\wedge\baeta] ) \\
  &=3\int_{M} \varphi^*(2[\theta\wedge \mu\wedge\nu\wedge\eta\wedge
  \btheta\wedge \bmu\wedge\bnu\wedge\baeta]) \\
 &=6\int_M [\theta\wedge \mu\wedge\nu\wedge\eta\wedge
  \btheta\wedge \bmu\wedge\bnu\wedge\baeta] \neq 0\, .
 \end{align*}
By Lemma \ref{lem:nonformal}, $\widetilde{M}$ is non-formal.
\end{proof}

\begin{remark}
The symplectic manifold $(\widetilde M,\widetilde \omega)$ is not
hard-Lefschetz. The $\ZZ_3$-invariant form $\nu\wedge \bnu$ on $M$
is not exact, but $\omega^2\wedge\nu\wedge\bnu =
2d(\theta\wedge\bmu\wedge\baeta\wedge\eta\wedge\bnu)$. This form
descends to the quotient $\widehat M$ and can be extended to
$\widetilde M$ via the process done at the end of the proof of the
previous theorem. Therefore the map $[\omega]^2 \colon
H^2(\widetilde M) \to H^6(\widetilde M)$ is not injective.

Cavalcanti \cite{Cav2} gave the first examples of  simply
connected compact symplectic manifolds of dimension $\geq 10$
which are hard Lefschetz and non-formal. Yet examples of
non-formal simply connected compact symplectic $8$-manifolds
satisfying the hard Lefschetz property have not been constructed.
\end{remark}

\bigskip

\noindent {\bf Acknowledgments.} We are very grateful to the referee
for useful comments that helped to simplify 
the exposition in Section \ref{sec:formality}. We also thank Dominic
Joyce for suggesting us to look at \cite{Gu} and Gil Cavalcanti,
Ignasi Mundet and John Oprea for conversations and helpful
suggestions. This work has been partially supported through grants
MCyT (Spain) MTM2004-07090-C03-01, MTM2005-08757-C04-02 and Project
UPV 00127.310-E-15909/2004.

{\small

\vspace{0.15cm}

\noindent{\sf M. Fern\'andez:} Departamento de Matem\'aticas,
Facultad de Ciencia y Tecnolog\'{\i}a, Universidad del Pa\'{\i}s
Vasco, Apartado 644, 48080 Bilbao, Spain.

 {\sl E-mail:} marisa.fernandez@ehu.es

\vspace{0.15cm}

\noindent{\sf V. Mu\~noz:} Departamento de Matem\'aticas, Consejo
Superior de Investigaciones Cient{\'\i}ficas, C/ Serrano 113bis, 28006
Madrid, Spain.

Facultad de Matem{\'a}ticas, Universidad Complutense de Madrid, Plaza
de Ciencias 3, 28040 Madrid, Spain.

 {\sl E-mail:} vicente.munoz@imaff.cfmac.csic.es

\end{document}